\documentstyle{amsppt}

\magnification=1200 \NoBlackBoxes \hsize=11.5cm \vsize=18.0cm
\def\sch{\text{sch}}

\def\inv{^{-1}}
\def\?{{\bf{??}}}

\def\A{\Bbb A}

\def\C{\Bbb C}

\def\P{\Bbb P}

\def\[{\big[}
\def\]{\big]}

\def\O{\Cal O}

\def\y{\bar{y}}

\def\Sym{\text{Sym}}

\def\rk{\text{rk}}

\def\m{\frak m}

\def\1/2{\frac{1}{2}}

\def\I{\Cal I}

\def\im{\text{im}}

\def\2{{[2]}}

%\nologo
\topmatter
\title Geometry on Nodal Curves
\endtitle
\author
Ziv Ran
%\thanks{\raggedright{
%Partially supported by NSA Grant MDA904-02-1-0094} }
\endauthor

\date October 12, 2002\enddate

\address University of California, Riverside\endaddress
\email ziv\@math.ucr.edu\endemail
%\urladdr http://www.math.ucr.edu/~ziv/papers/
%semicurv.pdf\endurladddr
\rightheadtext { Geometry on Nodal Curves}
\leftheadtext{Z. Ran}
\abstract
Given a family $X/B$ of nodal curves we construct
canonically and compatibly with base-change, via
an explicit blow-up of the Cartesian product $X^r/B$,
a family $W^r(X/B)$ parametrizing length-$r$ subschemes
of fibres of $X/B$ (plus some additional data). Though
$W^r(X/B)$ is singular, the important sheaves on it
are locally free, which allows us to study intersection
theory on it and deduce enumerative applications,
including some relative multiple point formulae,
enumerating
the
length-$r$ schemes contained simultaneously in some
fibre of $X/B$ and some fibre of a given map from $X$
to a smooth variety.

%We develop enumerative methods
%for the study of groups of points on
%a family of semistable curves
\endabstract
 \thanks \raggedright {Updates and corrections
available at $\underline{math.ucr.edu/\tilde{\ }
ziv/papers/geonodal.pdf}$}\linebreak
Research Partially supported by
NSA Grant MDA904-02-1-0094; reproduction and
distribution of reprints by US governement permitted.
\endthanks

\endtopmatter\document
One of the important facts which make geometry,
in particular enumerative geometry, on a smooth curve
$X$ relatively simply is the existence of simple
and quite tractable parameter spaces for
subschemes of $X$ of given length $r$,
be it the symmetric product $\Sym^r(X)$, which in fact
is isomorphic to the Hilbert scheme Hilb$^r(X)$,
or the Cartesian product $X^r$ which parametrizes
subshemes in an $(r!:1)$ fashion, and is
sufficient for many applications, especially
enumerative ones. {\it{mutatis mutandis}}, the same holds
for families of smooth curves. In \cite{R1},
the author studied
enumerative projective geometry for families of smooth
curves obtaining, {\it{inter alia,}}
a general {\it{relative
multiple point formula}}, i.e. a formula enumerating
the length-$r$ subschemes of the fibres of a given family
$X/B$ whose image under a given map
$$f:X\to Y$$
is a single reduced point.\par
Surprisingly, it seems these
ideas and results have yet to be extended to singular
curves and families of such. This paper
is a step in that direction in the case of nodal curves,
i.e. curves with only ordinary
double points as singularities.
 To a
family
$$\pi:X\to B$$
of nodal curves and a natural number $r$
 we shall associate canonically a family
$$\pi_r:W^r(X/B)\to B$$ which is functorial
in $B$, i.e. its formation commutes with base-change;
in fact, $W^r(X/B)$ is a canonical
 and explicit blowup of the Cartesian
fibre product $X^r/B$ (more explicitly, of
$W^{r-1}(X/B)\times
_BX$) in a suitable sheaf of ideals. We will show
that $W^r(X/B)$ admits a morphism to the relative Hilbert
scheme Hilb$^r(X/B)$. In fact,
we will subsequently show that
$W^r(X/B)$ is isomorphic to the relative
{\it{ ordered flag-
Hilbert scheme}}, which parametrizes chains
$$z_1\subset z_2...\subset z_r$$
where each $z_i$ is a length-$i$ subscheme of a fibre
of $X/B$, together with compatible orderings of the support
of each $z_i$.\par
The fact that $W^r(X/B)$ admits a morphism to Hilb
implies that it carries 'tautological' bundles (also
called secant bundles) $S^r(L)$, for any vector
bundle $L$ on $X$,  whose fibre at a point
is the restriction of $L$ on the corresponding scheme.
We will see that $S^r(L)$ can be analyzed conveniently
with exact sequences. This fact, together with
the fact that certain 'diagonal' divisors become
Cartier on $W^r(X/B)$ enables us to do some
intersection theory on these spaces and apply it
to enumerative questions (although the complete
intersection theory of the $W^r(X/B)$ is yet to be
worked out).\par
The remainder of the paper is largely
devoted to multiple-point
formulae. In the case of a map $f$ to a Grassmannian,
one can
define and enumerate a multiple-point scheme
$M_r(f)$ as a
more-or-less direct consequence of the existence of a
tautological bundle on $W^r(X/B)$ (whose fibre at a point
is the space of functions, or sections
of a vector bundle, the corresponding subscheme), and
its relation to the tautological bundles
on the Grassmannian. In the
case of a map to a general (smooth) variety $Y$,
this approach is not available
so instead we develop an essentially
iterative approach, viewing $M_r(f)$
essentially as a sublocus of (a blowup of)
$M_{r-1}(f)\times_BX$, we have
a natural mapping
$$M_{r-1}(f)\times_BX\to Y\times Y$$
and we can define $M_r(f)$ as the locus in the pullback of
the diagonal in $Y\times Y$ residual to the 'diagonal'
where the $r$th point coincides with one of the first
$r-1$. The fundamental class of $M_r(f)$ can then be
computed with the Fulton-MacPherson residual-intersection
formula [F].
\par
The paper is organized as follows. In a preliminary \S 1
we study finite schemes supported at a (1-dimensional)
node or in a neighborhood of a node in a family of curves;
this is basically a matter of elementary algebra.
In \S 2 we give the construction of the parameter spaces
$W^r(X/B)$, proceeding by induction. We then derive our
multiple-point formula for maps to Grassmannians and give
some applications, mainly for maps to $\P^2.$
In \S 3 we give the general multiple-point formula, allowing
maps to an arbitrary smooth target variety.\par
This paper is a continuation of \cite{R6} where more
particular results were obtained, essentially by
an {\it{ad hoc}} version of the methods of this paper.

\heading 1. The Hilbert scheme at a node
\endheading
In this section we will first study the punctual
Hilbert scheme of
length-$r$ schemes supported at the origin on a
germ of a node
$xy=0.$ Then we will study the full Hilbert
scheme of this germ.
Finally we will study the relative Hilbert scheme
of a family of
curves in the neighborhood a node. These results will form
the local foundation on which we shall
in \S 2 construct our global
parameter spaces $W^r(X/B)$.\par

We denote by $R$ the localization of the ring
$$\C[x,y]/(xy)=
\{a+\sum\limits_{i\geq 1}b_ix^i+c_iy^i\}$$
at its maximal ideal $(x,y)$.
Thus the formal completion
$$\hat{R}=\C[[x,y]]/(xy)$$
is isomorphic to the formal completion of the local
ring at any 1-dimensional ordinary node. We seek to
determine the Hilbert scheme
Hilb$^0_m(R)$ of colength-$m$
ideals in $R$ which, as is well known, is naturally
isomorphic to Hilb$^0_m(\hat{R})$.
\proclaim {Proposition 1.1}
(i)Any ideal $I<R$ of colength
$m$  is of one of the following types:
$$ {(c^m_i)}\ \ \ I^m_i(a)=
(y^{i}+ax^{m-i}), a\neq 0, i=1,...,m-1;$$
$$(q^m_i)\hskip2cm \ \ \ Q^m_i=(x^{m-i+1},y^{i}),
i=1,...,m.$$
(ii) The closure $C^m_i$ in the Hilbert scheme of
the set of ideals of type $(c^m_i)$ is a $\P^1$
and consists of the ideals of types $(c^m_i)$ or
$(q^m_i)$ or $(q^m_{i+1})$. In fact, we have
$$\lim\limits_{a\to 0}I^m_i(a)=Q^m_i,$$
$$\lim\limits_{a\to \infty}I^m_i(a)=Q^m_{i+1}.$$
(iii) The punctual Hilbert scheme Hilb$_m^0(R)$
is a rational chain
$$C^m_1\cup_{Q^m_2}C^m_2\cup\cdots\cup C^m_{m-1},$$
with nodes at $Q^m_2,...,Q^m_{m-1}$
and smooth elsewhere.
\linebreak
(iv) The only colength-$(m-1)$ ideal containing
$I^m_i(a)$ is $Q^{m-1}_i$; the only colength-$(m-1)$
ideals containing $Q^m_i$ are the $I^{m-1}_{i-1}(a)$
for $a\neq 0$ and their limits $Q^{m-1}_i$ and
$Q^{m-1}_{i-1}$.

\endproclaim
\demo{proof} Note that any nonzero nonunit $z\in R$ is
associate to a uniquely determined element of the form
$x^\alpha$ or $y^\beta$ or $x^\alpha+ay^\beta, a\neq 0,
\alpha,\beta>0,$
in which case we will say that $z$
is of type $(\alpha,0)$
or $(0,\beta)$ or $(\alpha,\beta)$, respectively. Note
also that for any ideal $I$ of colength $m$ we have
$$x^m, y^m\in I.$$
Now given $I$ of colength $m$, pick $z\in I$ of
minimal type $(\alpha,\beta)$,
with respect to the natural
partial ordering on types. Suppose to begin with that
$\alpha,\beta>0.$ Then note that
$$x^{\alpha+1}, y^{\beta+1}\in I,$$
and consequently $(\alpha,\beta)$ is unique:
indeed if $(\alpha',\beta')$ is also minimal
then we may assume $\alpha'>\alpha$, hence
$x^{\alpha'}\in I$, hence $y^{\beta'}\in I,$
contradicting minimality.
Hence $(\alpha,\beta)$ is unique
and it is then easy to see that the element
$z'=x^\alpha+ay^\beta\in I$ is unique as well,
so clearly $z'$ generates $I$ and $I$ is of type
$(c^m_\beta)$.\par
Thus we may assume that any minimal element of $I$
is of type $(\alpha,0)$ or $(0,\beta)$. Since
$x^m,y^m\in I$, $I$ clearly contains minimal
elements of type $(\alpha,0)$ and $(0,\beta)$,
and then it is easy to see that $I$ is of type
$(q^m_\beta)$. This proves assertion (i).
Since $I^m_i(a)$ contains $y^{i+1}, x^{m-i+1},$
assertion (ii) is easy, as is (iv).\par
For (iii), we note that the punctual Hilbert
scheme parametrizes punctual deformations of an
ideal $I<R$, parametrized by a local augmented artin
$\C-$algebra $S$. This means ideals
$$I_S<R_S:=S[x,y]_{(x,y)}$$
such that $R_S/I_S$ is $S$-free of rank $m$,
and which are punctual in the sense
that $I_S$ is contained in a unique ideal
$J$ that is maximal subject to
$$J\cap S=0$$ (in which
case
$J=(x,y)$ automatically).
The assertion is that this scheme is 1-dimensional
and smooth at $I^m_i(a)$ and nodal at $Q^m_i.$
 We consider the case of $Q^m_i$
as that of $I^m_i(a)$ is similar but simpler. Then
$I_S$ is generated by
$$f=x^{m+1-i}+f_1(x)+f_2(y),$$
$$g=y^i+g_1(x)+g_2(y)$$
where the $f_i,g_j$ are polynomials without constant
term (by punctuality) and with coefficients in $\m_S$.
If the low term of $f_1$ is $bx^j, j<m+1-i$ then
$$I_S\subseteq (x^{m+1-i-j}+b,y)$$
contradicting maximality. Therefore $j\geq m+1-i$
and changing by units we may assume
$$f=x^{m+1-i}+by^r$$
$$g=y^i+cx^s.$$
Note
%$$x^{m+2-i}, y^{i+1}\in I_S,$$
if $r\geq i$ then replacing $f$ by $f-bx^{r-i}g$ yields
a polynomial in $x$ only which by the above remarks
may be assumed to be $x^{m+1-i}$, so we can take
$b=0$; hence we may assume $r<i$ and similarly
$s<m+1-i$. On the other hand if $r<i-1$ then
$by^{i-1}\in I_S$ whereas by flatness,
$$1, x,...,x^{m-i}, y,...,y^{i-1}$$
must be an $S$- free basis of $R_S/I_S$.
Therefore $r=i-1, s=m-i$. Finally, since
$$cx^{m+1-i}+bcy^{i-1}, by^i+bcx^{m-i}\in I_S,$$
subtracting these from $g,f$ respectively
we get a contradiction to punctuality unless $bc=0.$
Thus, punctual $S-$deformations of $Q^m_i$ are
parametrized by a pair $b,c,\in\m_S$ with $bc=0,$
i.e. by the local $\C-$algebra homomorphisms
of the local ring of a node to $S$, so the punctual
Hilbert scheme itself is a 1-dimensional node locally
at $Q^m_i$, as claimed.
\qed\enddemo
Next we determine the structure of the full Hilbert
scheme of $R$:
\proclaim{Proposition 1.2} The Hilbert scheme
$\text{Hilb}_m(R)$  is a chain
$$D^m_0\cup D^m_1\cdots D^m_{m-1}\cup D^m_m$$
where each $D^m_i$ is a smooth and $m-$dimensional germ
supported on $C^m_i$ for $i=1,...,m-1$ or $Q^m_i$
for $i=0,m$; for $i=1,...,m-1, D^m_i$ meets its
neighbors $D^m_{i\pm 1}$ transversely
in dimension $m-1$ and meets no other
$D^m_i.$\endproclaim
\demo{proof} Clearly Hilb$_m(R)$ is a germ supported on
Hilb$^0_m(R)$, so  this is a a matter of determining the
scheme structure of Hilb$_m(R)$ at each point of
Hilb$^0_m(R)$, which may be done formally by testing on
Artin local algebras. Again, we shall do so
at $Q^m_i$ as the case of $I^m_i(a)$ is similar and
simpler. Given $S$ artinian local augmented, a  flat
$S$-deformation of $I=Q^m_i$ is given by an ideal
$$I_S=(f,g),$$
$$f=x^{m+1-i}+f_1(x)+f_2(y),\tag 1.1$$
$$g=y^i+g_1(x)+g_2(y),$$
where $f_i,g_j$ are polynomials with coefficients
in $\m_S$, and such that $R_S/I_S$ is $S-$ free
of rank $m$, in which case it is clear by Nakayama's
Lemma that
$$1, x,..., x^{m-i}, y ,..., y^{i-1}$$
is a free basis for $R_S/I_S$.
It is easy to see that we may assume
$f_1, g_1$ are of degree $\leq m-i$ and $f_2,g_2$
are of degree $< i$ and $f_2,g_2$ have no constant
term.  Let's write
$$f_1(x)=\sum\limits_0^{m-i} a_jx^j, f_2(y)=\sum\limits
_1^{i-1} b_jy^j,\tag 1.2$$
$$g_1(x)=\sum\limits_0^{m-i} c_jx^j, g_2(y)=\sum\limits
_1^{i-1} d_jy^j.\tag 1.3$$
Now obviously
$$yf-b_{i-1}g\equiv 0\equiv xg-c_{m-i}f\mod I_S.$$ Writing
these elements out in terms of
$1,x,...,x^{m-i},
y,...,y^{i-1}$
yields relations among
$1,x,...,x^{m-i},
y,...,y^{i-1}$. Since the latter elements form an $S$-free
basis of $R_S/I_S$, those relations must be trivial,
in other words we have exact equalities rather than
congruences:
$$yf-b_{i-1}g= 0= xg-c_{m-i}f.$$
Writing this out yields the following  identities
$$b_j=b_{i-1}d_{j+1}, j=1,...,i-2,$$
$$b_{i-1}d_1=a_0,$$
$$b_{i-1}c_j=0, j=0,...,m-i,\tag 1.4$$
$$c_j=c_{m-i}a_{j+1}, j=0,...,m-i-1,$$
$$c_{m-i}a_0=0,$$
$$c_{m-i}b_j=0, j=1,...,i-1.$$
Conversely, suppose the
relations (1.4) are satisfied, i.e.
$$yf-b_{i-1}g= 0= xg-c_{m-i}f.\tag 1.4'$$
By Nakayama's Lemma, $1, x,...,x^{m-i}, y,...,y^{i-1}$
generate $R_S/I_S$, hence to
show
(1.1) defines a flat family it suffices
to show these elements
admit no nontrivial $S$-relations mod $I_S$.
To this end,, suppose
$$u_{m-i}(x)+v_{i-1}(y)=A(x,y)f+B(x,y)g\tag 1.5$$
where $u(x),...$ are all polynomials in the indicated
variables and of the indicated degees
(if any) with coefficients in $S$ and $v$ has no constant
term (in fact it clear that then $A,B$ have coefficients
in $\m_S$). Then the relations (1.4')
 allow us
to rewrite (1.5) as
$$u_{m-i}(x)+v_{i-1}(y)=A'(x)f+B'(y)g,\tag 1.6$$
and comparing coefficients in (1.6) we see
directly that
$$A'=B'=u_{m-i}=v_{i-1}=0,$$
hence there are no nontrivial relations, as claimed.\par
Thus the Hilbert scheme is embedded in the
space of the variables $$a_1,...,a_{m-i}, d_1,...,d_{i-1},
b_{i-1},c_{m-i},$$ i.e. $\A^{m+1}$, and defined by the
relation
$$b_{i-1}c_{m-i}=0.\tag 1.7$$
Thus it is a union of 2 smooth $m-$dimensional components
meeting transversely in a smooth $(m-1)-$dimensional
subvariety.
\par
In the case $I=I^m_i(a)=(y^i+ax^{m-i})$, a similar
analysis shows that an $S$-deformation of $I$ is
given by
$$I_S=(y^i+\tilde{a}x^{m-i}+f_1(x)+f_2(y))$$
where
$$\tilde{a}\in S, \tilde{a}\equiv a\mod\m_S,$$
$$f_1(x)=\sum\limits_0^{m-i-1} a_jx^j, f_2(y)=\sum\limits
_1^{i-1} b_jy^j,$$
$$a_j, b_j\in\m_S,$$
and via $(\tilde{a}, a_0,...,a_{m-i-1}, b_1,...,b_{i-1})$
we may identify Hilbert scheme locally with $\A^m.$
\qed
\enddemo
Finally we consider the relative local situation,
i.e. that of a germ of a (1-parameter) family of
curves with smooth  total space specializing
to a node. Thus set
$$\tilde{R}= \C[x,y]_{(x,y)},    B=\C[t]/_{(t)},
$$
and view $\tilde{(R)}$ as a $B-$module via $xy=t$. As is
well known, this is the versal deformation of the node
singularity $xy=0$, so any family of nodal curves is
locally a pullback of this.
\proclaim{Proposition 1.3} The relative Hilbert scheme
$\text{Hilb}^m(\tilde{R}/B)$ is smooth, formally
$(m+1)-$dimensional
as a total space.\endproclaim
\demo{proof} The relative Hilbert scheme parametrizes
length-$m$ schemes contained in fibres of
Spec$\tilde{(R)}\to$
Spec$(B)$. This means ideals $I_S<\tilde{R}_S$ of
colength $m$ containing $xy-s$ for some $s\in\m_S$,
such that $\tilde{R}_S/I_S$ is $S-$free. The analysis
of these is virtually identical to that contained in
the proof of Proposition 1.2, except that the relation
$b_{i-1}c_{m-i}=0$ gets replaced by
$$b_{i-1}c_{m-i}=s.$$
This means precisely that the relative Hilbert scheme
is the subscheme of the affine space of the variables
$a_1,...,a_{m-i},d_1,...,d_{i-1}, b_{i-1}, c_{m-i}, t$
defined by the relation
$$b_{i-1}c_{m-i}=t,\tag 1.8$$
hence is smooth as claimed.\qed\enddemo
The local analysis immediately yields
some conclusions for the Hilbert scheme of
a nodal curve:
\proclaim{Corollary 1.4} Let $C_0$ be a curve
with only $k$ nodes as singularities and $c$ irreducible
components. Then\par\noindent
(i)\ \  $\text{Hilb}^m(C_0)$ is reduced and
has precisely
$\binom{m+c-1}{m}$ components, the general element
of each of which corresponds to a reduced subscheme
of the smooth part of $C_0;$
\par\noindent (ii)\ \ let $I$
be a point of
$\text{Hilb}^m(C_0)$ having colength
$m_i$ at the $i$-th node of $C_0;$
then locally at $I$,

$\text{Hilb}^m(C_0)$ is a cartesian product of
$k$ factors, each of which is a 2-component normal
crossing of dimension $m_i, i=1,...,k$,
or a point if $m_i=0$, times
a smooth factor.
\par\noindent (iii)\ \ the fibre of the cycle map
$$cyc:\text{Hilb}^m(C_0)\to \Sym^m(C_0)$$
over a cycle having multiplicity $m_i$ at the
$i$th node is a product of
1-dimensional rational chains of length $m_i-1.$
\endproclaim
\demo{proof} It is clear from the explicit analysis
in the proof of Proposition 1.2 that any subscheme
of $C_0$ deforms to a reduced subscheme supported on
the smooth part. Such subschemes are
parametrized by an open dense subset of the
symmetric product $\Sym^m(C_0).$ This clearly yields
(i) and (ii), while (iii) follows from the fact that
the fibres of $cyc$ are products of punctual
Hilbert schemes.\qed\enddemo
\heading 2. Parameter spaces\endheading
Let
$$\pi:X\to B$$
be a flat family of curves over an irreducible
variety, with all fibres nodal and generic fibre
smooth. Our purpose here is to construct a natural
and explicit birational modification $W^r(X/B)$
of the relative
cartesian product $X^r/B$, which will serve as our
basic 'configuration space' on which to do
enumerative geometry.\par
We begin by considering in explicit detail
the cases of small $r$; subsequently the construction
wil proceed by induction. Obviously, set
$$W^1(X/B)=X/B;$$ we next consider
the case $r=2$. Note that, locally at a critical
point $p$ of $\pi$ (i.e. a singular point of a fibre
$X_b=\pi\inv(b)$), our family is analytically
equivalent to a subvariety of $\A^2\times B$
given by
$$xy=a\tag 2.1$$
where $x,y$ are coordinates on $\A^2$ and $a$ is
a function on $B$, which may be also be viewed
as a mapping to the base ($=\A^1$) of the versal
deformation of a node, pulling back  the germ at $p$
of $X/B$. Note that $X$ is smooth at $p$ iff
$da_p\neq 0$. If $\dim B=1,\ a$ may be
taken to be of the form
$a=t^k$, where $t$ is a local coordinate on $B$ and
$k\geq 1$ and $k=1$ iff the
$X$ is smooth at $p$.\par
Now consider the fibre square $X^2/B$ and let
$D\subset X^2/B$ be the diagonal. Note that $D$ is a
Cartier divisor {\it{except}} at points $(p,p)$
where $p$ is a fibre singularity; there, $X^2/B$ is
given locally by
$$x_1y_1=x_2y_2=a,\tag 2.2$$
and $D$ is given by the 2 equations
$$x_1=y_1, x_2=y_2;\tag 2.3$$
moreover if the total space $X$ is smooth,
then $X^2/B$ is smooth (over $\C$) except at
those points $(p,p)$. We let
$$b_2:W^2(X/B)\to X^2/B$$ denote
the blowup of $D\subset X^2/B$, i.e. the blowup,
locally at all points $(p,p)$ as above of the ideal
(2.3), and let
$$\Delta^2\subset W^2(X/B)$$
be the exceptional divisor, i.e. $b_2\inv(D)$ with
the natural scheme structure endowed by the blowing-up
process (not to
be confused with the exceptional {\it{locus}}, i.e.
the locus of all points of $W^2(X/B)$ where $b_2$
is not an isomorphism locally. It will suffice to
analyze the situation locally along the exceptional
locus of $b_2.$\par
In terms of coordinates, over a neighborhood of each
$(p,p)$, $W^2(X/B)$ is covered by 2 open affines
denoted $U^2_{x,p}, U^2_{y,p}$
or just $U^2_x, U^2_y$ if $p$ is understood.
The coordinate ring of $U^2_{x,p}$
is generated over that of $X^2/B$ by a symbol
$[(y_2-y_1)//(x_2-x_1)]$ subject to the evident relation
$$(x_2-x_1)[(y_2-y_1)//(x_2-x_1)]=y_2-y_1.$$
Note that
$$x_1(x_2-x_1)[(y_2-y_1)//(x_2-x_1)]=x_1(y_2-y_1)$$
$$=y_2(x_1-x_2),$$
hence we may, and shall, write $[(y_2-y_1)//(x_2-x_1)]$
as $-[y_2//x_1]$; similarly, the same may also
be written as $-[y_1//x_2],$ therefore when the index
range is understood we may write the same as
$-[y//x]$. Similar comments
apply over the other open $U^2_{y,p}$ where a generator
$[x_1//y_2]=[x_2//y_2]= [x//y]$ is defined and of course,\
on the overlap $U^2_{x,p}\cap U^2_{y,p}$ we have
$$[y_2//x_1][x_1//y_2]=1.$$
Henceforth we shall denote $U^2_{x,p}, U^2_{y,p}$
respectively as $U([y//x]), U([x//y])$.\par
Thus the exceptional locus of $b_2$ consists of
a $\P^1$ over each point $(p,p)$ as above.
Moreover it is easy to see that if $X$ is smooth,
then so is $W^2(X/B)$; in fact, if $X$ is a smooth
surface then $(p,p)$ is just a 3-fold ordinary
double point and $b_2$ is one of its 2 small resolutions.
It is useful to note that this resolution may be obtained
determinantally: i.e. via (2.2) we obtain
locally a mapping
$$M: X^2/B\to M^1_{2\times 2},$$
$$M=\left [\matrix x_1&x_2\\y_2&y_1\endmatrix\right],$$
to the space of $2\times 2$ matrices of rank
$\leq 1$
and $W^2(X/B)$ is just obtained by taking fibre product
of $X^2/B$ via $M$ with the canonical determinantal
resolution of $M^1_{2\times 2}$,
$$R^1_{2\times 2}=\{ (A,B)\in M^1_{2\times 2}\times\P^1:
BA=0\}.$$
Next, we claim that the natural rational map
of $X^2/B$ to the Hilbert scheme Hilb$_2(X/B)$,
assigning a pair $p\neq q$ to the ideal
$\Cal I_{p,q}<\O_X$ lifts to a morphism
$$J_2:W^2(X/B)\to \text{Hilb}_2(X/B).$$
Clearly it suffices to check this locally along
the exceptional locus. To this end it suffices
to note that in the open subset
$U^2_{x,p},$ $J_2$ is given explicitly by sending a point
with coordinates $(x_1,x_2,[y_2//x_1])$ to the ideal
$$J_2(x_1,x_2,[y_2//x_1])=$$
$$ (xy-a, (x-x_1)(x-x_2),
y-y_1+[y_2//x_1](x-x_1))\tag 2.4$$
where, of course, we set
$$y_1=x_2[y_1//x_2], y_2=x_1[y_2//x_1],
a=x_1y_1=x_2y_2=x_1x_2[y_2//x_1].$$
Indeed, it is elementary that the formula (2.4)
defines an ideal of colength
2 whose cosupport contains
$(x_1,y_1), (x_2,y_2)$, which implies our assertion.
The case of $U^2_{y,p}$ is similar.\par
Next, before tackling the case of general $r$
we will study in detail the case $r=3.$
Let
$$\Gamma^3\subset W^2\times_BX$$
denote the pullback of the tautological subscheme
of Hilb$_2(X/B)$ via the map $J_2$. Note that
$$\Gamma^3=\Gamma^3_1\cup\Gamma^3_2\tag 2.5$$
where $$\Gamma^3_i\subset W^2\times_BX$$ is the graph
of the natural projection $W^2\to X$: indeed this
follows easily from the fact that the LHS of (2.5)
is flat over $W^2$, the RHS is reduced and both sides
agree generically over $W^2.$ In particular, we see that
$\Gamma^3$ is reduced. We define
$$W^3=W^3(X/B):=B_{\Gamma^3}(W^2(X/B)\times_BX),$$
i.e. the blowup of $W^2\times_BX$ in $\Gamma^3,$
with natural map
$$b_3:W^3\to W^2\times_BX.$$
Let $\Delta^3\subset W^3$ be the natural exceptional
(Cartier) divisor supported on $b_3\inv(\Gamma^3)$
and $\Delta^3_i=b_3\inv(\Gamma^3_i), i=1,2$ its
components (as Weil, in fact nonCartier divisors).
It is easy to see that all these divisors are reduced,
as is $W^3.$\par
To analyze this construction, we work over
$U([y//x]).$
There, note that the expression $(y_3-y_1)(x_3-x_1)$
viewed, e.g. as a function on
$$V:=U([y//x])\times_BU,$$
is divisible by $x_2$: indeed setting
$$R=y_2-[y//x]x_3-y_3+y_1,$$
it is easy to check that $x_2R=(y_3-y_1)(x_3-x_1).$
I claim next that the ideal of $\Gamma^3$ in $V$ is
generated by $(x_3-x_1)(x_3-x_2), R$: indeed the subscheme
$\Gamma'$
defined by the latter projects isomorphically to the
subscheme of the '$x$-axis' defined by $(x_3-x_1)(x_3-x_2)$,
hence is flat over $W^2$, and of course $\Gamma^3$
is also flat over $W^2.$ Since $\Gamma^3$ and $\Gamma'$
clearly coincide generically over $W^2,$ they coincide, as
claimed.\par
Thus we see that $b_3\inv(V)$ may be realized in the standard
way as a subscheme of $V\times\P^1$ and as such is covered
by the two standard opens pulled back from $\P^1$.
One of these is the domain of regularity of
 the rational function
$$-[y//x^2]:=\frac{R}{(x_3-x_1)(x_3-x_2)},$$
and it is easy to check that as rational functions,
$$[y//x^2]=\frac{y_i}{x_jx_k}$$
whenever $\{ i,j,k\}=\{ 1,2,3\}$, justifying the notation;
we denote this open by $U([y//x^2])$. Note, trivially,
that regularity of $[y//x^2]$ implies that of $[y_i//x_j],
\forall i,j,$ so $U([y//x^2])$ indeed lies over
$U([x//y])$. On the other standard open of $V\times\P^1$
the function
$$[x^2//y]=\frac{x_jx_k}{y_i}$$
is regular, however its domain of regularity does
not lie entirely over $U([y//x])$ (nor, for that matter,
entirely over $U([x//y]).$\par
Analogous comments apply to the part of $W^3$ over
$U([x//y])$ which gives rise to symbols $[x//y^2],
[y^2//x]$. Note by construction that the regularity
domains
$$U([y//x^2]), U([y^2//x], [x^2//y]):
=U([y^2//x])\cap U( [x^2//y]), U([x//y^2])$$
already cover $W^3(X/B)$.\par
I claim next that the natural rational map $J_3$
of $W^3(X/B)$ to the relative Hilbert scheme
Hilb$_3(X/B)$ is a morphism. This is a local assertion,
and is moreover obvious except at points of $W^3$
lying over $(p,p,p)\in X^3/B$, so it suffices to check
it on opens
$$U([y//x^2]), U([y^2//x], x^2//y]),
U([x//y^2])$$ as above. Over $U([y//x^2]),$
it is easy to see that the expression
$$(y-y_1)(x-x_2)
(x-x_3),$$ considered as a function on $W^3\times_BX$,
i.e. taken modulo $xy-x_1y_1,$ is divisible by
$x_2x_3$: explicitly, setting
$$R=([y_2//x_3]-[y//x^2])x^2+([y_1//x_3]+[y_1//x_2])x
+y-(y_1+y_2+y_3)$$
(recall that $[y_2//x_3]=x_1[y//x^2]$
 etc.), we have
$$x_2x_3R\equiv (y-y_1)(x-x_2)(x-x_3)\mod xy-x_1y_1.$$
Then on $U([y//x^2]), J_3$ takes a point with coordinates
$x_i,y_j$ to the ideal
$$J_3=((x-x_1)(x-x_2)(x-x_3), R, xy-x_1y_1).$$
Since the latter ideal evidently has colength 3 over $W^3$,
this makes $J_3$ a morphism over $U([y//x^2]).$
The case of $U([x//y^2])$ is similar. Over
$U([y^2//x], [x^2/y]),$ it is elementary
to check as above that, always modulo $xy-x_1y_1,$
$(y-y_1)(y-y_2)(x-x_3)$ is divisible by $x_3$ and
$(x-x_1)(x-x_2)(y-y_3)$ is divisible by $y_3$ and the
ideal
$$J_3=(\frac{(y-y_1)(y-y_2)(x-x_3)}{x_3},
\frac{(x-x_1)(x-x_2)(y-y_3)}{y_3}, xy-x_1y_1)$$
has colength 3 over $W^3$ and yields the map to
Hilb (cf. \S 1).\par
For general $r, W^r(X/B)$ is naturally constructed
by induction. It is convenient to summarize the
construction and its main properties as follows.
\proclaim{Theorem 2.1} Define
$$b_r: W^r(X/B)\to W^{r-1}(X/B)\times_BX$$
inductively
as the blowup of
 in the
canonical subscheme
$$\Gamma^r\subset W^{r-1}(X/B)\times_BX$$
corresponding
to the morphism
$$J_{r-1}:W^{r-1}(X/B)\to \text{Hilb}_{r-1}(X/B),$$
and let
$$w_r:W^r(X/B)\to X^r/B$$
be the natural map.
Then\par
(i) $W^r(X/B)$ is reduced, and is irreducible if $X$ is;
$\Gamma^r$ is reduced and has $r-1$ irreducible components
each isomorphic to $W^{r-1}$; as functor
of the family $X/B$, $W^r(X/B)$
commutes with base-change;
\par
(ii) the natural rational map
$$J_r:W^r(X/B)\to{\text{Hilb}}_r(X/B)$$
is a morphism;
\par
(iii) for each node
$p$ of $X/B$, a neighborhood $N$ of $w_r\inv(p,...,p)$
in $W^r=W^r(X/B)$ carries rational functions
$$[y^i//x^{r-i}]=[y^I//x^J], [x^i//y^{r-i}]=[x^I//y^J]$$
where $I\subset [1,r]$ is any index-set of cardinality
$i$ and complement $J$ and the domains of regularity
$$U([y//x^{r-1}]),...,U([y^i//x^{r-i}],
[x^{r-i+1}//y^{i-1}]),
...,U([x//y^{r-1}])$$
form a covering of $N;$\par
(iv) setting
$$P^r_i=\prod\limits_{j=1}^i(y-y_j)\prod\limits_{j=i+1}
^r (x-x_j)\in \O_{W^r}[x,y]/(xy-x_1y_1), \tag 2.6$$
over $U([y^i//x^{r-i}]),\ P^r_i$
is divisible by $x_{i+1}\cdots x_r$ and
over $U([x^{r-i}//y^i]),\ P^r_i$
is divisible by $y_{1}\cdots y_i$; over
$U([y^i//x^{r-i}],[x^{r-i+1}//y^{i-1}]),$
the map $J_r$ is given by
$$J_r=(xy-x_1y_1, \frac{P^r_i}{x_{i+1}\cdots x_r},
\frac{P^r_{i-1}}{y_1\cdots y_{i-1}});$$
\par (v) $W^{r+1}$ is covered by open sets over each of
which either\par
(a)$[y^i//x^{r+1-i}]$ is regular and
$\frac{P^r_{i-1}(x_{r+1}, y_{r+1})}{y_1\cdots y_{i-1}}$
is an equation
for $\Delta^{r+1}$ (called an '$x-$based equation); or\par
(b)$[x^{r-i+1}//y^{i}])$ is regular and
$\frac{P^r_i(x_{r+1},y_{r+1})}{x_{i+1}\cdots x_r}$
is an equation
for $\Delta^{r+1}$ (called a '$y-$based equation).
\endproclaim
\demo{proof} For $r\leq 3$ all of this has already
been proven, so we may assume it holds for $r-1.$
First (i) is clear from the fact that $W^r$ is an
iterated blowup of the cartesian product $X^r/B$, while
$\Gamma^r$ coincides with the union of the graphs
(over $B$) of
the coordinate projections $p_i:W^{r-1}\to X$
(proof as in the $r=3$ case). For the rest, we may
as before work over
$$U([y^i//x^{r-1-i}], [x^{r-i}//y^{i-1}])\subset
W^{r-1}$$
so in a suitable open set in $W^{r-1}\times_BX$,
the ideal of $\Gamma^r$ is generated by
$$\frac{P^{r-1}_i}{x_{i+1}\cdots x_{r-1-i}},
\frac{P^{r-1}_{i-1}}
{y_1\cdots y_{i-1}}$$
(where we plug in $(x_r,y_r)$ for $(x,y)$).
Thus the blowup (in the part under consideration) is
covered by two standard opens. In the first we have
the regular function
$$\frac{P^{r-1}_i}{x_{i+1}\cdots x_{r-1-i}}/
\frac{P^{r-1}_{i-1}}
{y_1\cdots y_{i-1}}=\frac{y_r-y_i}{x_r-x_i}
\frac{y_1\cdots y_{i-1}}{x_i\cdots x_{r-1-i}}$$
and it is easy to see as before that this coincides as
rational function with
$$-\frac{y_1\cdots y_i}{x_{i+1}\cdots x_r}$$
and, for that matter, with any $-y^I/x^J$ as in (iii)
so we may denote it by $$-[y^i//x^{r-i}].$$
It is also easy to see as before that this standard open
coincides with the regularity domain of this function
so we denote it by $U([y^i//x^{r-i}])$.
Similarly we get a rational function $[x^{r-i}//y^i]$.
\par
Now we can prove (iii). Given $z\in W^r$,
such that $w_r(z)$ is close to $(p,...,p)$,
 we may assume $z$
projects to
$$z'\in U([y^i//x^{r-1-i}], [x^{r-i}//y^{i-1}])\subset
W^{r-1},$$
and in particular
$$\frac{y_1\cdots y_iy_r}{x_{i+1}\cdots x_{r-1}}=
[y^{i+1}//x^{r-1-i}],$$$$
\frac{x_i\cdots x_{r-1}x_r}{y_1\cdots y_{i-1}}=
[x^{r-i+1}//y^{i-1}]$$
is regular at $z$.
As we have seen, either $$[y^i//x^{r-i}]\ {\text{ or}}\
[x^{r-i}//y^i]$$ are regular at $z$. Therefore either
$$z\in U([y^i//x^{r-i}],[x^{r-i+1}//y^{i-1}])$$
or
$$z\in U([x^{r-i}//y^i],[y^{i+1}//x^{r-1-i}]).$$
Thus (iii) is proved.\par
We will now prove (iv). To begin with,
note that the question is local (on $W^r$, {\it{
a fortiori}} on $X/B$) so we may assume $X/B$ is the versal
family $xy=t$ over $B=\A^1$ (actually we
just need that $B$ is integral). Now
 we first show that over
$U([y^i//x^{r-i}]),\ P^r_i$
is divisible by $x_{i+1}\cdots x_r$; the proof that
over $U([x^{r-i}//y^i]),\ P^r_i$
is divisible by $y_{1}\cdots y_i$
is similar. Now
 $P^r_i$ (a polynomial in $x,y$
 subject to the relation
 $xy=x_1y_1$) is a sum of terms of the form
$$M_{i-a}(y)M_{r-i-b}(x)y^{a}x^{b}$$
where the $M$'s are monomials in distinct variables
$y_1,...,y_i,x_{i+1},..., x_r$ of the indicated degrees.
If $a\leq b$, use the relations
$$xy=x_jy_j, \forall j$$
to rewrite this term as
$$x^{b-a}y_1\cdots y_iM_{r-i-b+a}(x),$$
which is clearly divisible as claimed. If $b\leq a$
this term can be rewritten as
$$y^{a-b}M_{i-a+b}(y)x_{i+1}\cdots x_r,$$
which is even more clearly divisible as claimed.\par
Note that the above calculation shows
 $P^r_i/x_{i+1}\cdots x_r$ and
$P^r_{i-1}/y_1\cdots y_{i-1})$ can be written respectively
as
$$y^i+f_1(x)+f_2(y),$$
$$x^{r+1-i}+g_1(x)+g_2(y)$$
where $f_1,f_2$ have degree $\leq r-i$ and $g_1,g_2$
have degree $\leq i-1$ and all have regular
functions as coefficients. By the proof of Proposition 1.2,
we see that for $J_r$ as defined to yield a morphism
$$J_r:W^r(X/B)\to \text{Hilb}^r(X/B)$$
is equivalent to
 certain identities (1.4) among the coefficients of
the $f_i$ and $g_j$.
Since $J_r$ clearly coincides {\it{generically}}
with the evident rational map,
these identities hold generically, hence they hold period,
so $J_r$ is indeed a lifting of the evident rational
map as morphism to Hilb.
This completes the proof of (iv). Finally,
in light of the fact that
$$\frac{P^r_{i-1}(x_{r+1}, y_{r+1})}{y_1\cdots y_{i-1}}
/ \frac{P^r_i(x_{r+1},y_{r+1})}{x_{i+1}\cdots x_r}
=[x^{r-i+1}//y^{i}],$$
(v) follows directly from (iv) and
the definition of blowup
(the two opens in question are the standard affine
opens of $\P^1$ over
$U([y^i//x^{r-i}],[x^{r-i+1}//y^{i-1}])$.
This completes the proof
of Theorem 2.1.
\enddemo
A posteriori, we can identify $W^r(X/B)$ as the
'ordered flag-Hilbert' scheme: recall that the flag
Hilbert scheme fHilb$_m(X/B)$ parametrizes $\O_B-$
chains of ideals
$$ I_m<...<I_1<\O_X$$
where $\O_X/I_j$ is $\O_B-$locally free of length $j$.
We set
$$\text{ofHilb}_m(X/B)=\text{fHilb}(X/B)\times_BX^m/B.$$
Note that the various maps $J_j$ together yield a morphism
$$\eta_r: W^r(X/B)\to\text{ofHilb}_r(X/B).$$
\proclaim{Corollary 2.2} The map $\eta_r$ is an
isomorphism.\endproclaim
\demo{proof} As usual, it suffices to prove that the
fibre of $\eta_r$ over an ordered
 flag supported at a point $p$ that is a relative node
 reduces (scheme-theoretically) to a point.
 By induction, it then suffices to prove
 that a fibre of the natural map
$$\zeta_r:
W^r\to W^{r-1}\times_BX\times_B\text{Hilb}_r(X/B)$$
is a point.
Note that by construction, if $w\in W^r$ lies over
$$w'\in U([y^i//x^{r-1-i}], [x^{r-i}//y^{i-1}])$$
then a fibre of the map
$$b_r: W^r\to W^{r-1}\times_BX$$ already
is coordinatized locally by either
$$[y^i//x^{r-i}]\ {\text{or}}\
[x^{r-i}//y^i].$$
Now by Proposition 1.2 and the computations in
the proof if Theorem 2.1(iv) is is clear
that the corresponding point in Hilb$_r$ is either
$$I^r_i(a),\ \text{with}\ \ a=[y^i//x^{r-i}]$$
if $a\neq 0$ (note that $a$ is the coefficient of
$x^{r-i}$ in $P^r_i/x_{i+1}\cdots x_r$), or $Q^r_i$ if
$$[y^i//x^{r-i}]=[x^{r-i+1}//y^{i-1}]=0.$$
In either case, the computations show that
$[y^i//x^{r-i}]$ or $[x^{r-i}//y^i]$ can be read off
from Hilb$_r$, consequently $\zeta_r$ is an isomorphism,
hence so is $\eta_r$.\qed

\enddemo
Unfortunately, $W^r$ is never smooth if $r\geq 3$,
as soon as $X/B$ has singular fibres,
but we still have
\proclaim{Proposition 2.3} If $B$ is smooth then
$W^r(X/B)$ is normal and Cohen-Macaulay.\endproclaim
\demo{proof} We first show inductively that $W^r$ is
$R_1.$ For $r=1$ this is clear. Inductively, if
$W^{r-1}$ is $R_1$ then clearly so is $W^{r-1}\times_BX$.
Moreover the blowup map
$$b_r:W^r\to W^{r-1}\times_BX$$
has at most $\P^1$ fibres and those only over a
codimension-3 locus. Hence $W^r$ is $R_1.$\par
For CM, we again argue inductively. If $W^{r-1}$
is CM, then since $X$ is CM and $W^{r-1}\times_BX$
is a locally complete intersection in
$W^{r-1}\times X$, it is CM as well.
Finally $W^r$ is locally a hypersurface in
$W^{r-1}\times_BX\times\P^1$, it is CM as well.\qed
\enddemo
Now let
$$\Gamma^{r+1}\subset W^r\times_BX$$
be the 'tautological divisor' (=pullback of universal
divisor over the Hilbert scheme via $J_r$). For any
locally free sheaf $L$ on $X$, set
$$S^r(L)=p_{W^r*}(p_X^*(L)\otimes\O_{\Gamma^{r+1}}),
\tag 2.7$$
which is clearly locally free of rank $r.\rk(L)$,
and which we call
the $r-$th {\it{ secant bundle}} associated to $L$.
It was introduced in the smooth case by
Schwarzenberger \cite{Sc}.
Note that $W^r(X/B)$ is not symmetric with respect
to permuting the factors, but still there are
projections 'to the first $s$ factors', for all
$s\leq r$:
$$\gamma^{r,s}:W^r\to W^{s}.$$
We also set
$$\gamma^r=\gamma^{r,r-1}.$$ Also, denote by
$$\Delta^r\subset W^r$$
be the exceptional divisor of $b_r$, i.e.
$b_r^*(\Gamma^r)$, which is by definition of
blowup a Cartier divisor. As we saw before in the
$r=3$ case, we have for any $r$ that $\Gamma^r$
splits up as
$$\Gamma^r=\bigcup\limits_{i=1}^{r-1}\Gamma^r_i$$
with each $\Gamma^r_i$, being the graph of the $i-$th
coordinate projection to $X$, is isomorphic to
$W^{r-1}$ and in particular is reduced always
and irreducible if $X$ is. Similarly $\Delta^r$
splits
$$\Delta^r=\bigcup\limits_{i-1}^{r-1}\Delta^r_i$$
with each $\Delta^r_i$ reduced and birational
to $\Gamma^r_i$ (and in general a non-Cartier
divisor on $W^r$).
\proclaim{Proposition 2.4}
There is an exact sequence of vector bundles
on $W^r(X/B)$
$$0\to p_r^*(L)\otimes\O(-\Delta^r)\to S^r(L)\to
\gamma^{r*}(S^{r-1}(L))\to 0 \tag 2.8$$
\endproclaim
\demo{proof} There is clearly a natural surjection
$$S^r(L)\to\gamma^{r*}(S^{r-1}(L))$$
whose kernel $K$ is locally free of rank $\rk(L)$
and moreover admits a generically injective map
$$k:K\to p_r^*(L).$$
Clearly, $k$ vanishes at each generic point of
$\Delta^r$, hence factors through
$p_r^*(L)\otimes\O(-\Delta^r)$. It is easy to see that
the factored map
$$K\to p_r^*(L)\otimes\O(-\Delta^r)$$
is an isomorphism in codimension 1 of line bundles,
hence an isomorphism since $W^r$ is pure-dimensional.
\qed\enddemo
Thus $S^r(L)$ has a natural filtration with quotients
$$(\gamma{r,s})^*(p_s^*(L)\otimes\O(-\Delta^s)),$$
and consequently we have
\proclaim{Corollary 2.5} The total Chern class
$c(S^r(L))$ satisfies
$$c(S^r(L))=\prod\limits_{i=1}^r
(1+L_i-(\gamma^{r,i})^*\Delta^i)\tag 2.9$$
where $L_i$ denotes (the class of) $p_i^*(L)$
and $\Delta^1=0$.\endproclaim
This result can be used to give a {\it{ multiple-point
formula}} for maps to a Grassmannian. Let
$X/B$ be as above and
$$f:X\to G$$
be a map to a Grassmannian $G=G(k,m+k)$, with tautological
sub- and quotient bundles $S_G, Q_G$ respectively. Set
$$S_X=f^*S_G, Q_X=f^*Q_G.$$
Note that on $W^r(X/B)\times G$ we have a natural
map
$$p_2^*(S_G)\to p_1^*S^r(Q_X)$$ which vanishes at a point
$(w,g)$ iff $f$ maps the scheme corresponding to $w$
to the reduced point scheme $g$. We call the latter
locus $M^+_{r}(f).$ If we define a bundle
$S^{r,1}(L)$, for any $L$, by the exact sequence
$$0\to S^{r,1}(L)\to S^r(L)\to (\gamma^{r,1})^*(L)\to 0$$
then clearly the zero-scheme of the induced map
$$p_2^*(S_G)\to p_1^*(\gamma^{r,1})^*(Q_X)$$
is just the graph of the composite
$$f\circ p_1:W^r\to G.$$ Therefore $M^+_{r}(f)$ projects
isomorphically to its image $M_r(f)$ in $W^r$ and
$M_r(f)$ is the zero-scheme of an analogous map
$$p_1^*(S_X)\to S^{r,1}(Q_X).$$

Then we have shown
\proclaim {Theorem 2.6} In the above situation,
$M^+_{r}(f)$ and $M_r(f)$ have a natural scheme structure as
zero-scheme of vector bundle maps and when they have their
expected codimension, viz. $rkm$ or $(r-1)km$,
the cohomology class

$$[M_r(f)]=
c_{(r-1)km}(p_1^*(S_X^*)\otimes S^{r,1}(Q_X)).\tag 2.10$$

\endproclaim
In case $G=\P^m=G(m,m+1),$ the formula (2.10) simplifies
somewhat. Let's write
$$L_i=(f\circ p_i)^*(\O_G(1), M_1=(f\circ p_1)^*(S_X).$$
Then we get

$$[M_r(f)]=
\prod\limits_{i=2}^r
c_m(M_1^*\otimes(L_i-(\gamma^{r,i})^*(\Delta^i)))$$
$$=\prod\limits_{i=2}^r(\sum\limits_{j=0}^m
L_1^j(L_i-(\gamma^{r,i})^*(\Delta^i))^{m-j}).\tag 2.11$$
Let's consider the case $m=2,r=3.$ Then the RHS of (2.11)
yields
$$(L_1^2+L_1(L_2-\Delta^2)+(L_2-\Delta^2)^2)
(L_1^2+L_1(L_3-\Delta^3)+(L_3-\Delta^3)^2).\tag 2.12$$
We want to compute the image of this on $X=W^1$.
To this end we must first compute
the image on the second factor
on $W^2$ via $\gamma^{3,2}$. This computation
follows formally from the following
formulae, where we set
$$f^i_0={\text{ fibre of}}\ W^i(X/B) {\text{
over}}\ 0\in B,$$$$\ d=\deg(f^1_0), b= \pi_*L^2 :$$\par
1. $\gamma^{32}_*(L_1^2)=0$ (obvious)\par
2. $\gamma^{32}_*(L_1L_3)=dL_1$\par
3. $\gamma^{32}_*(L_1\Delta^3)=2L_1\ $ (since $\Delta^3$
has
2 components mapping birationally to $W^2$)\par
4. $\gamma^{3,2}_*(L_3^2)=\pi_2^*(b)$ (obvious)\par
5.$ \gamma^{3,2}_*(L_3\Delta^3)=L_1+L_2$ (cf. (3))\par
6. $\gamma^{3,2}_*((\Delta^3)^2)=-K_1-K_2+2\Delta^2$
where $K_i=p_i^*(K_{X/B})$; to show this we may work
off the (codimension-2) exceptional locus of the natural
birational map
$$W^2\to X^2/B$$
and its inverse image in $W^3$; on this open set,
 $\Delta^3$ consists of
2 components $\Delta^3_1, \Delta^3_2$, each a pullback of
the diagonal via the $p_{13}, p_{23}$ projections,
which meet in a locus projecting isomorphically to
$\Delta^2\subset W^2$.\par
[Actually analogues of formulae (1-6) hold for any
$gamma^{r.r-1}, r\geq $ in place of $\gamma^{3,2}$, where the
analogue of (6) is
$$(\gamma^{r,r-1})_*((\Delta^r)^2)=-\sum\limits_{i=1}^r
K_i+2\sum\limits_{i=2}^{r-1}(\gamma^{r-1,i})^*(\Delta^i)
$$
where we set $\gamma^{r-1,r-1}=$identity.]\par
The result is
$$\gamma^{3,2}_*((L_1^2+L_1(L_3-\Delta^3)+(L_3-\Delta^3)^2))$$
$$=(d-4)L_1-2L_2-K_1-K_2+2\Delta^2+\pi_2^*(b)\tag 2.13$$
Thus the image of (2.12) via $\gamma^{3,2}$ is
$$(L_1^2+L_1(L_2-\Delta^2)+(L_2-\Delta^2)^2)\cdot$$$$
((d-4)L_1-2L_2-K_1-K_2+2\Delta^2+\pi_2^*(b)) \tag 2.14 $$
The computation of this product is straightforward,
as is that of its image in $X$, noting that the normal
bundle to $\Delta^2$ in $W^2$ is just
$-b_2^*K_1+E$ where
$$b_2:W^2\to X^2/B$$
is the natural blowup map and $E$ its exceptional locus,
which is a divisor on $\Delta^2$; indeed the restriction
of $b_2$ over the diagonal
$$\Delta_X=X\subset X^2/B$$
is just the blowup of the critical locus of $\pi$, which we
denote by $\sigma$
(cf. \cite{R6}), therefore
$$\gamma^{2,1}_* (\Delta^2)^2= -K,$$
$$\gamma^{2,1}_* (\Delta^2)^3=
\gamma^{2,1}_* (K-E)^2= K^2-\sigma .$$

Notice that the image on $X$ of $M_3(f)$ is geometrically
twice the locus of points contained in a relative triple
point of $f$, while the image on $B$ of the same is 6
times the locus of fibres containing a triple point.
Thus writing out the product yields
\proclaim{Theorem 2.7} Let $\pi:X\to B$
be a family of connected
nodal
curves of arithmetic genus $g$
and $f:X\to\P^2$ a morphism. Then the virtual class
on $X$ of the locus of points contained in a relative
triple point of $f$ is
$$N_{3,X}(f)=\frac{1}{2}((3d^2-18d+24+6g-6)L^2
+(18-3d)KL+4K^2-2\sigma);\tag 2.15$$
the locus in $B$ of fibres containing a relative triple
point of $f$ has virtual class
$$N_{3,B}(f)=\frac{1}{3} \pi_*(N_{3,X}(f)).$$
Here $L=f^*\O(1), K=\omega_{X/B}, d=\deg f({\text{fibre}}),
\sigma=$ critical locus of $f$.\endproclaim
Let us finally specialize to the case where
$X/B$ is normalization of the pencil of rational
curves in $\P^2$ through $3d-2$ assigned generic points.
In this case all the ingredients of   (2.15) have
been computed recursively
before (cf. \cite{R4, R5} and references therein):
$L^2=N_d$, the number of rational curves of degree
$d$ through $3d-1$ generic points,
$$K=-2s_1-m_1f^1_0+R_1$$
where $s_1$ is a section of $X/B$ contracted by $f$ to
a point, $R_1$ is the sum of all fibre components
disjoint from $s_1$, and $m_1$ is a certain
recursively computed number. Hence
$$L.K=-dm_1+L.R_1, K^2=-\sigma$$
and $\sigma$ coincides with the number of reducible
fibres of $X/B$, also recursively computed.
We conclude
\proclaim{Corollary 2.8} The number of rational curves
of degree $d$ in $\P^2$ having a triple point and going
through $3d-2$ generic points is (cf.\cite{R6})
$$N_{d,3}=
\frac{1}{2}((d^2-6d+10)+(d-6)(dm_1-L.R_1)-3\sigma)\tag 2.16$$
\endproclaim
\remark{Example 2.9} $N_{4,3}=60$, a number first computed
by Zeuthen and rederived with modern methods by
Kleiman and Piene
\cite{KP} (I am grateful to Steve Kleiman for
this reference). See \cite{KP,R6} for some similar
examples. \endremark
When $f$ is a map to $\P^m$, one is interested classically
not only in the relative multiple-point loci of $f$
but also in its relative {\it{multisecant loci}}, that is
the locus of length-$r$ subschemes of fibres
whose $f-$image is contained in a linear $\P^k$.
This locus can be enumerated by the above results,
as the $r-$fold locus of the natural projection
$$I_X\to \Bbb G(k,m):=G$$
where $I_X$ is the incidence variety, i.e.
$$I_X=\{(x,L):f(x)\in L\}\subset X\times G.$$
But it is simpler and more direct to enumerate this
locus as follows.
Set
$${\text{Sec}}^r_k(f)=\{(z,L):f(\sch(z))\subset L\ \
{\text{as schemes}}\}\subset W^r(X/B)\times G\tag 2.17$$
where $\sch(z)$ is the subscheme of $X$ corresponding
to $z$. Clearly ${\text{Sec}}^r_k(f)$ is just the zero-
scheme of the natural map
$$p_2^*(Q_G^*)\to S^r(L)$$
where $Q_G$ is the tautological
quotient bundle (of rank $m-k$)
on $G$ and $L=f^*(\O(1))$. Thus we conclude
\proclaim{Theorem 2.9} For a family of nodal curves
$X/B$ and a morphism
$$f:X\to\P^m,$$ the virtual locus
on $W^r(X/B)\times G$
of relatively
$r-$secant $k$-planes to $X/B$ in $\P^m$  is given by
$$[{\text{Sec}}^r_k(f)]=c_{r(m-k)}
(p_1^*(S^r(L))\otimes p_2^*(Q_G)).\tag 2.16$$
If the RHS of (2.16) is nonzero (resp. noneffective),
then the locus of relative $r-$secant $k-$planes is
nonempty (resp. of dimension larger than the expected,
viz. $\dim(B)+r+(k+r+2)(m-k)$).\qed
\endproclaim
Note that the projection of
${\text{Sec}}^r_k(f)$ to $W^r(X/B)$ coincides
with the locus where the natural map
$$H^0(\O_{\P^m}(1))\otimes\O_{W^r(X/B)}\to S^r(L)$$
has rank at most $k+1$, and consequently can be enumerated
directly via Porteous' formula [F].

\heading 3. Multiple-point Formulae\endheading
The purpose of this section is to state and prove
a general {\it{relative multiple point formula}} for families
of nodal curves $X/B$, i.e. a formula enumerating the
length-$r$ subschemes in fibres of $X/B$ mapped to a point
by a given map
$$f:X\to Y$$
to a general smooth variety.
See \cite{F, F1, K} and references therein for
other work on multiple-point formulae. \par
In \2 we saw that in case $Y$ was a Grassmannian we could
write down, and enumerate, a multiple-point scheme
as a zero-scheme of a suitable vector bundle on $W^r(X/B)$.
This was a more-or-less direct consequence
of, on the one hand the standard representation
of the diagonal
of $Y$ as a zero-set of a vector bundle and,
on the other hand,
the existence of a tautological (secant)
bundle on $W^r$. For general $Y$ such a representation
of the diagonal is not available, and consequently
we will adopt an inductive, or recursive, approach,
defining the $r$-th multiple-point scheme in terms of
the $(r-1)$st.\par
The first step is to define
inductively the appropriate {\it{multiple point schemes}}
$$M_r\subseteq W^r(X/B).$$
Thus set
$$M_1=W^1(X/B)=X.$$
The inductive step is provided by the following
\proclaim{Lemma 3.1} Locally over $(\gamma^{r+1})\inv (M_r)$,
$f(xy_1)-f(xy_{r+1})$ is divisible by an equation
of $\Delta^{r+1}$.\endproclaim
The proof of the Lemma is based on a kind of
'finite-difference calculus' that we now explain.
Let $f=f(x,y)$ be a vector-valued 'function'
(which may be regular, rational, a (formal or convergent)
power series, $C^\infty$, etc.). Define operators
$D_{x_{21}}, D_{y_{21}}$ by
$$D_{x_{21}}(f)(x_1,y_1,x_2,y_2)=\frac{f(x_1,y_1)-
f(x_2,y_1)}{x_1-x_2},$$
$$D_{y_{21}}(f)(x_1,y_1,x_2,y_2)=\frac{f(x_2,y_2)-
f(x_2,y_1)}{y_2-y_1}$$
where the functions on the RHS are 'regular' in the same
sense as $f$.
These may be iterated, leading to operators
$D_{*_{mi}\cdots *_{21}}$ with each $*=x$ or $y$,
the value being a function of $m$ pairs $(x_i,y_i)$
which may be independent or subjected to to me relations-
in the case we need
the relations
$$x_1y_1=x_2y_2=\cdots=x_my_m$$
will be imposed. For convenience we may write
$$D_{x_{321}}=D_{x_{32}}D_{x_{21}},$$
etc.
The basic properties of these operators to used are the
following\par
(1) {\it{commutativity}}, e.g.
$$D_{y_{31}}D_{x_{21}}=D_{x_{23}}D_{y_{31}};$$
etc.\par
(2){\it{product rule}}, e.g.
$$D_{x_{21}}(fg)(x_1,y_1,x_2,y_2)=$$
$$D_{x_{21}}f.g(x_1,y_1)+f(x_2,y_2)D_{x_{21}}g;$$
(3) {\it{chain rule}}, e.g.
$$\frac{f(x_1,y_1)-f(x_2,y_2)}{x_1-x_2}=
D_{x_{21}}f+\frac{y_1-y_2}{x_1-x_2}D_{y_{21}}f$$
with $$\frac{y_1-y_2}{x_1-x_2}$$ only a rational function;
if the relation $x_1y_1=x_2y_2$ holds, the latter can
be rewritten as
$$D_{x_{21}}f-[y//x]D_{y_{21}}f$$
where, as before
$$[y//x]=y_1/x_2=y_2/x_1.$$

\demo{proof of Lemma}
As usual, it suffices to prove this
in a neighborhood of a point lying over $(p,...,p)$
for a fibre node $p$. We will first prove it
in a neighborhood of a
point in $U([y_{r+1}//x_1\cdots x_r]).$ As we have seen,
in such a neighborhood a local equation for
$\Delta^{r+1}$ is given by
$(x_1-x_{r+1})\cdots(x_r-x_{r+1})$
so me must show
$$\frac{f(xy_1)-f(xy_{r+1})}
{(x_1-x_{r+1})\cdots(x_r-x_{r+1})}$$
is regular.
In fact, we
will prove more precisely that in such a neighborhood
we have
$$\frac{f(xy_1)-f(xy_{r+1})}{(x_1-x_2)\cdots (x_r-x_{r+1})}
=(-1)^{r+1}D_{x_{(r+1)r...21}:y_1}f+R_{1...r(r+1)}\tag  3.1 $$
where $R_{1...r(r+1)}$ is a sum of terms of the form
$$\pm [y_{r+1}//x^{I_1}]\cdots [y_{r+1}//x^{I_k}]
D_{x_*y_*x_*...:*}f.$$ For $r=1$ this is easy,
with the usual freshman calculus technique of
'subtracting and adding $f(x_2,y_1):$
$$\frac{f(xy_1)-f(xy_{2})}{(x_1-x_2)}=
D_{x_{21}:y_1}f+\frac{y_1-y_2}{x_1-x_2}
D_{y_{21}:x_2}f$$
$$=D_{x_{21}:y_1}f-[y_2//x_1]
D_{y_{21}:x_2}f.$$
So assume ? is true with $r+1$ replaced by $r$,
so that on $M_r$ we have
$$(-1)^{r}D_{x_{(r)(r-1)...21}:y_1}f=
-R_{1...(r-1)r}\tag  3.2 $$

Then we also have an analogous
relation to (3.2) replacing $r$ by $r+1$:
$$\frac{f(xy_1)-f(xy_{r+1})}
{(x_1-x_{r+1})\cdots (x_{r-1}-x_{r+1})}
=(-1)^{r}D_{x_{(r+1)(r-1)...21}:y_1}f
+R_{1...(r-1)(r+1)}\tag  3.3 $$
(a priori the coefficients in (3.3) are only rational
but, being all regular multiples of
$[y_{r+1}//x_1\cdots x_r]$, are regular in the open
where considered). Then from (3.1) we get
$$\frac{f(xy_1)-f(xy_{r+1})}{(x_1-x_2)\cdots (x_r-x_{r+1})}=$$

$$\frac{1}{x_r-x_{r+1}}((-1)^{r}D_{x_{(r+1)(r-1)...21}:y_1}f
+R_{1...(r-1)(r+1)})=\tag  3.4 $$
$$\frac{1}{x_r-x_{r+1}}((-1)^{r}D_{x_{(r+1)(r-1)...21}:y_1}f
-(-1)^{r}D_{x_{(r)(r-1)...21}:y_1}f+$$$$
(-1)^{r}D_{x_{(r)(r-1)...21}:y_1}f
+R_{1...(r-1)(r+1)})=\ $$
$$\frac{1}{x_r-x_{r+1}}((-1)^{r}D_{x_{(r+1)(r-1)...21}:y_1}f
-(-1)^{r}D_{x_{(r)(r-1)...21}:y_1}f$$$$
-R_{1...(r-1)r}
+R_{1...(r-1)(r+1)})\tag 3.5$$
(the latter equality by (3.2)).
Now clearly,
$$\frac{1}{x_r-x_{r+1}}((-1)^{r}D_{x_{(r+1)(r-1)...21}:y_1}f
-(-1)^{r}D_{x_{(r)(r-1)...21}:y_1}f
$$$$=(-1)^{r+1}D_{x_{(r+1)r(r-1)...21}:y_1}f$$
and by 'derivative of a product' formulae as in elementary
calculus, we see that
$$\frac{1}{x_r-x_{r+1}}(-R_{1...(r-1)r}
+R_{1...(r-1)(r+1)})$$
has the form $R_{1...(r+1)}$ as claimed.\par
Finally, in the general case we use Theorem 2.1 (v) and
symmetry which allow us to work in an open set where
$[x^{r-i}//y^i]$ and $[y^{i+1}//x^{r-i}]$
are regular and
$$E=\frac{P^r_i(x_{r+1},y_{r+1})}{y_{1}\cdots y_i}$$
is an equation for $\Delta^{r+1}.$
We recall that here and in similar situations,
$[x^{r-i}//y^i]$ is short for
$(\gamma^{r+1})^*([x^{r-i}//y^i])$, i.e.
$x_{i+1}\cdots x_r/y_1\cdots y_i$.
Working inductively,
starting with
$$(f(x_1,y_1)-f(x_{r+1},y_{r+1}))y_1/
(y_1-y_{r+1}),$$ we can check that
$$(f(x_1,y_1)-f(x_{r+1},y_{r+1}))/E$$
 is a linear combination
of expressions $D_{*\cdots *}f$ with
each $*=x_{ab}$ or $y_{cd}$ and
coefficients $C$ which
are products of $x_i,y_j$'s,
$[x^{r-i}//y^i]$ and $[y^{i+1}//x^{r-i}]$,
hence are regular (in fact, counting each $x,y$
index in $D_{*\cdots *}$ as (-1) , the total bidegree
if $C.D_{*\cdots *}$ is $(-(r-i),0)).$
This completes the proof of
Lemma 3.1.\enddemo
Now the Lemma gives us an inclusion
$$i:(\gamma^{r+1})^*(\I_{M_r})\subseteq \O_{W^{r+1}}(
-\Delta^{r+1})$$
and we define a subscheme $M_{r+1}$ by
$$\I_{M_{r+1}}=\im (i\otimes\O(\Delta^{r+1}))\subset
\O_{W^{r+1}}.$$
This is what's usually
called the {\it{residual scheme}}
 to $\Delta^{r+1}$
in $(\gamma^{r+1})^*({M_r}).$\par
To state our multiple-point formula we need some
further notation. Let $X/B$ be as above and
$$f:X\to Y$$
be a morphism to a smooth $m$-dimensional variety.
For $k\geq 2$ set
$$\mu_k(f)= ((f\circ p_1)\times (f\circ p_k))^*(\Delta_Y)
-\Delta^k\{\frac{(f\circ p_1)^*(c(T_Y)}{1+\Delta^k}\}_{m-1}
\tag 3.6.$$
Also let
$$m_r(f)=\prod\limits_{k=2}^r (\gamma^{r,k})^*(\mu_k(f))
\tag 3.7$$
where
$$\gamma^{r,k}:W^r(X/B)\to W^k(X/B)$$
is the natural map
\proclaim{Theorem 3.2} Let $X/B$ be a family of nodal
curves and $$f:X\to Y$$ a morphism to a smooth
$m$-dimensional
variety.
Then\par
(i) there is a natural scheme structure $M_r(f)$ on
the locus in $W^r(X/B)$ of points whose associated
scheme is mapped by $f$ to a reduced point;\par
(ii) if $B$ is irreducible,
$M_r(f)$ is locally defined
by $(r-1)m$ equations on $W^r(X/B)$ hence is
purely at least $(\dim(B)+r-(r-1)m)-$
dimensional and if equality holds then
$$[M_r(f)]=m_r(f).
\tag 3.8 $$\endproclaim
\demo{proof} The proof is by induction on $r$, and
clearly holds for $r=1$, so assume the result holds for
$r=1$. For (i) it would suffice that a point
$z_r\in W^r$ is in the support of
$M_r(f)$ (as constructed above) iff the
scheme corresponding to $z_r,$ i.e.
$$\sch(z_r):=\O_X/J_r(z_r)$$
is mapped by $f$ to a reduced point. Clearly, it
suffices to prove this when $\sch(z_r)$ is punctual,
supported at a fibre node $p$. Choosing coordinates
on $Y$, we may represent $f$ as a vector-valued function
with $f(p)=0.$ Let $z_{r-1}$ be the projection of
$z_r$ to $W^{r-1}.$ By Proposition 1.1(iv), either
$\sch(z_r)$ is of type $I$ and $\sch(z_{r-1})$ is of
type $Q$ or vice versa. We assume the former as the latter
is similar and simpler. Thus we assume
$$J_r(z_r)=(y^j+ax^{r-j}),
J_{r-1}(z_{r-1})=(y^j,x^{r-j}).$$
By construction of $W^r$, this means that local
equations for the tautological subscheme $\Gamma^r$
over a neighborhood of $z_r$ are
$$P^r_{j-1}/y_1\cdots y_{j-1}, P^r_j/x_{j+1}\cdots x_r$$
while equations for $\Gamma^{r-1}$ over a neighborhood
of $z_{r-1}$ are
$$P^{r-1}_{j-1}/y_1\cdots y_{j-1},
P^{r-1}_j/x_{j+1}\cdots x_{r-1},$$
and that a local equation for $\Delta^r$ in a neighborhood
of $z_r$ is
$$E = P^{r-1}_{j-1}(x_r,y_r)/y_1\cdots y_{j-1}$$
%$$=\frac{P^{r-1}_{j-1}/y_1\cdots y_{j-1}}{
%P^{r-1}_j/x_{j+1}\cdots x_{r-1}};$$
moreover, $[y^j//x^{r-j}]$ is regular at $z_r$,
indeed
$$[y^j//x^r-j]=\frac{P^{r-1}_j/x_{j+1}\cdots x_{r-1}}
{P^{r-1}_{j-1}/y_1\cdots y_{j-1}}$$
and
$$a=[y^j//x^{r-j}](z_r)).$$
Now by induction $f$ maps $\sch(z_{r-1})$ to the origin,
which means that, as a function on $W^r$, we can write
$$f= bP^{r-1}_{j-1}/y_1\cdots y_{j-1}+
cP^{r-1}_j/x_{j+1}\cdots x_{r-1}.$$
Thus
$$f/E=b+c[y^j//x^r-j]$$
so
$$(f/E)(z_r)=b+ac.$$
Thus $(f/E)(z_r)=0$ iff $f$ vanishes on $\sch(z_r)$,
completing the proof of (i).\par
In (ii), the assertion about the number of equations is
clear from the definition. As for the assertion about the
cohomology class, it is clear from the Fulton-MacPherson
residual- intersection formula [F] provided both $M_r$
and $M_{r-1}$ have their expected dimensions. In the
general case, let $C_1,...,C_k$ be the irreducible
components of $M_{r-1}$. By Fulton's theory, the is
a cycle $U_i$ of dimension $\dim(B)+r-1-(r-2)m$
on each $C_i$ such that
$$\sum [U_i]=m_{r-1}.$$
Since $M_r$ is locally defined by $m$ equations over
$(\gamma^r)\inv(M_r)$ but still has its expected dimension,
it follows that the contribution of each oversize
component $C_i$ to $M_r$ is empty, and in particular
$$(\gamma^r)^*(U_i).\mu_r(f)=0.$$
So these oversize components contribute nothing to either
$M_r$ or $m_r$, so \?? still holds.\qed

\enddemo
\proclaim{Corollary 3.3} In the situation of Theorem
3.2, if $m_r(f)\neq 0$ then $M_r(f)$ is nonempty.
\endproclaim
Let us say that a smooth $m$-dimensional variety $Y$
is {\it{pseudo-
Grassmannian}} with bundle $G$ if the diagonal
$$\Delta_Y\subset Y\times Y$$
is a zero-scheme of $G$.
\remark {Examples} (i)Clearly a Grassmannian
has this property, with
$$G=p_1^*S^*\otimes Q$$
where $S$ and $Q$ are respectively the tautological
sub- and quotient bundles. \par
(ii) Trivially, any curve is pseudo-Grassmannian\par
(iii) Generally, a product of pseudo-Grassmannians is
pseudo-Grassmannian, therefore any product of curves and
Grassmannians is pseudo-Grassmannian.\endremark
Then the proof
of Theorem 3.2 yields directly the following refinement.
\proclaim{Theorem 3.2 bis} In the situation of
Theorem 3.2, if moreover $Y$ is pseudo-Grassmannian
with bundle $G$ then $M_r(f)$ is a zero-scheme of
$$\bigoplus\limits_{i=2}^r
((f\circ p_1)\times (f\circ p_i))^*(G)(-(\gamma^{r,i})^*
\Delta^i)\tag 3.9$$
\endproclaim

Consider the case $r=2, m=3.$ Thus
we have a family of nodal curves mapping to $Y$
and are enumerating
the relative multiple points of their images in $Y$
(at least if we assume that a general fibre of $X/B$
is smooth and embeds in $Y$ and that every fibre maps in
with degree 1). Then \? yields (writing $f_i=f\circ p_i$):
$$[M_2(f)]=(f_1\times f_2)^*(\Delta_Y)
-((\Delta^2)^3-(\Delta^2)^2f_1^*K_Y+\Delta^2f_1^*c_2(Y)).
\tag 3.10$$
By the calculations in \cite{R6},
we have as in the proof of Corollary 2.8,
$$p_{1*}(\Delta^2)^3=K^2-\sigma,$$
$$p_{1*}((\Delta^2)^2f_1^*K_Y)=-K.f^*(K_Y),\
p_{1*}=f^*(c_2(Y))$$
where as before $K=\omega_{X/B}$ and $\sigma$
is the critical locus of $\pi$. Thus we obtain
\proclaim{Corollary 3.4} For a family of nodal curves
$X/B$ mapping via $f$ to a smooth 3-fold $Y$, the
virtual locus on $X$ of relative double points of $f$
is
$$[N_{2,X}(f)]=p_{1*}(f_1\times f_2)^*(\Delta_Y)-
(K^2-\sigma+K.f^*K_Y+f^*c_2(Y))\tag 3.11$$
\endproclaim
The expression $p_{1*}(f_1\times f_2)^*(\Delta_Y)$
(which is just a number if $B$ is 1-dimensional)
may be evaluated in various ways. For example, working
in singular cohomology over $\C$, let $(\alpha_i)$ be
a homogeneous basis for the total cohomology $H^*(B)$ and
let $(\alpha_i^*)$ be the dual basis. Then the class
of the diagonal $\Delta_B$ in $B\times B$ is given by
$$[\Delta_B]=\sum \alpha_i\otimes\alpha_i^*$$
setting
$$\beta_i=f_*(\pi^*(\alpha_i)),
\beta_i^*=f_*(\pi^*(\alpha_i^*))$$
we have
$$f_*(p_{1*}(f_1\times f_2)^*(\Delta_Y))=
\sum \beta_i.\beta_i^*\tag 3.12$$
(which coincides with $p_{1*}(f_1\times f_2)^*(\Delta_Y)$
or $f_1\times f_2)^*(\Delta_Y)$ when they are of top degree,
i.e. numbers. Note also that when $B$ is a curve, we have
$$[\Delta_B]\equiv [B]\otimes [pt]+[pt]\otimes B
\mod H^1\otimes H^1$$ so if $H^3(Y)=0$ then ?
reduces to $2f_*([f_0]).f_*([X]).$ Finally note that one
customarily denotes
$$\pi_*(K^2)=\kappa, \pi_*(\sigma)=\delta.$$
Thus we have
\proclaim {Corollary 3.5} In the above situation,
suppose\par
(i) a generic fibre of $X/B$ is smooth and embedded via
$f$;\par
(ii)$\dim(B)=1$;\par
(iii) $H^3(Y)=0.$\par
Then the virtual number of double of relative double
points of $f$ is given by
$$n_2(f)=\frac{1}{2}(2f_*([f_0]).f_*([X])
-\kappa+\delta-K.f^*K_Y
-f^*c_2(Y)).\tag 3.13$$
\par In particular, if $n_2(f)\neq 0$ then $f$ does not
embed all fibres of $X/B$ and if $n_2(f)<0$ then $f$
has degree $>1$ on some fibre.
\endproclaim
As a special case, we recover a result from \cite{R6}.
We use the notation developed there; in particular
$N_d^{\text{red}}(a_.)$ denotes the number of
{\it{reducible}} rational
curves of degree $d$ in $\P^3$ satisfying the incidence
conditions indicated by $(a.)$  and $m_1=-s_1^2$
where $s_1$ is the section of $X/B$ corresponding
to an incident linear subspace of codimension
$a_1$. Both these numbers are recursively
computable. See \cite{R4,R5} for more details.
\proclaim{Corollary 3.6} With the notations of \cite{R4,R5},
the number of singular rational curves of
degree $d$ through $a$ generic points and
$4d-2a-1$ generic lines in $\P^3$ is
$$(d-2)N_d(3^a2^{4d-2a})+N_d^{\text{red}}
(3^a2^{4d-2a-1})
-2dm_1+2L.R_1\tag 3.14$$
if $a>0$ and
$$(d-2)N_d(3^a2^{4d-2a})+N_d^{\text{red}}
(3^a2^{4d-2a-1})
-2dm_1+2L.R_1-4N_d(2^{4d-2a-2}3^{a+1})\tag 3.13$$
if $4d-2a-1>0.$
\endproclaim
\demo{proof} The RHS of \? yields, with
$(a.)=(3^a2^{4d-2a-1})$ or $(2^{4d-2a-1}3^a)$ that,
in the notation of \cite{R4},
$$n_2(f)=
\frac{1}{2}(2dN_d(3^a2^{4d-2a})-K^2+N_d^{\text{red}}(a.)
+4LK-4N_d(3^a2^{4d-2a})).$$
The intersection calculus developed in \cite {R4} shows that
$$LK=-2N_d(...,a_i+1,...)-dm_i+L.R_i$$
for any $i$ (we take $i=1$). This gives the result
(note that $N_d(4,..)=0$, whence the absence
of this term from (3.12).\enddemo
\remark {Remark 3.7} Theorem 3.2, as well as the other
multiple-point results in this paper, admit straightforward
generalizations to the relative case, where $Y$ is replaced by
a smooth morphism
$$\rho:Y\to B$$
and $f$ is a $B$-morphism, i.e. the following diagram
commutes
$$\matrix X&&\overset{f}\to{\to}&&Y\\
&\pi\searrow&&\swarrow\rho&&\\
&&B.&&
\endmatrix$$
Note that the 'absolute' case discussed above becomes
a special case of the relative case by replacing $Y$ by
$Y\times B\to B$. In the relative case
 The factors  $\mu_k(f)$ are replaced by
$$\mu_k(f/B)=(f_1\times_B f_k)^*(\Delta_{Y/B})-
\Delta^k\{\frac{f_1^*(c(T_{Y/B}))}{1+\Delta^k}\}_{m-1}
\tag 3.6'$$
where $\Delta_{Y/B}$ is the diagonal in
$Y\times_BY$ and $T_{Y/B}$ is the relative
or vertical tangent
bundle of $Y\to B$ (which coincides with the normal
bundle of $\Delta_{Y/B}$ in $Y\times_BY$) and $m$ is the
relative dimension of $Y/B$.
With (3.6') in place of (3.6), the analogue of Theorem
3.2 and its consequences hold. The analogous generalizations
of Theorem 2.6 and its consequences, which concern
maps to projective and Grassmannian bundles, also hold.
The proofs are the same, because for a $B-$map $f$,
multiple-point loci involve only the 'vertical' coordinates
of $Y$ over $B$.

\endremark

%\vfill\eject

\enddocument